\documentclass[11pt,a4paper]{amsart}

\usepackage{mathtools}
\usepackage{enumitem}
\setlist{itemsep=1mm}
\usepackage{hyperref}
\usepackage{url}
\usepackage[textwidth=27mm, textsize=scriptsize, color=green!20]{todonotes}

\usepackage{pdfsync}

\theoremstyle{definition}

\newtheorem{proposition-definition}[theorem]{Proposition-Definition}

\DeclareMathOperator{\diag}{diag}

 \newcommand{\iu}{{\mathrm{i}\mkern1mu}}
 \newcommand{\FS}{{\mathrm{FS}}}
\usepackage[OT2,OT1]{fontenc}

\makeatletter
\@namedef{subjclassname@2020}{\textup{2020} Mathematics Subject Classification}
\makeatother

\begin{document}

\title[Proskuryakov's determinant 308]{A note on the determinant 308 in
  Proskuryakov's linear algebra book}

\author[Di Scala]{Antonio J. Di Scala}
\address{Dipartimento di Scienze Matematiche, Politecnico di Torino. Corso Duca degli Abruzzi 24,
  10129 Torino, Italy}
\email{antonio.discala@polito.it}

\author[Sombra]{Mart{\'\i}n~Sombra}
\address{Instituci\'o Catalana de Recerca
  i Estudis Avan\c{c}ats (ICREA). Passeig Llu{\'\i}s Companys~23,
  08010 Barcelona, Spain  \vspace*{-2.5mm}}
\address{Departament de Matem\`atiques i
  Inform\`atica, Universitat de Barcelona. Gran Via 585, 08007
  Bar\-ce\-lo\-na, Spain}
\email{sombra@ub.edu}

\date{\today}   \thanks{Sombra was
  partially supported by the MINECO research project
  PID2019-104047GB-I0. Di Scala is member of CrypTO, GNSAGA of INdAM and of DISMA Dipartimento di Eccellenza MIUR 2018-2022.}

\begin{abstract}
  We put in evidence and correct a mistake in the formula for the
  determinant 308 in Proskuryakov's linear algebra book.  We apply
  this formula to reprove the well-known fact that the Fubini-Study
  metric on the complex projective space is Einstein.
\end{abstract}

\maketitle

This short note is motivated by a mistake in the formula for the
interesting determinant 308 in Proskuriakov's classical book of linear
algebra problems.  We checked several of its many editions including
the some of first ones and of the more recents~\cite{MR0219548, binom}
as well as the translations \cite{Reverte,MR797072}, and noticed that
the mistake has not been corrected.

Problem 308 asks to compute the determinant
\begin{equation}
  \label{eq:4}
  P308 = \det \begin{bmatrix}
x_1     & a_1 b_2 & a_1 b_3     & \cdots & a_1 b_n \\
a_2 b_1 & x_2     & a_2 b_3    & \cdots & a_2 b_n \\
a_3 b_1 & a_3 b_2 & x_3        & \cdots & a_3 b_n \\
\vdots  & \vdots  & \vdots     & \ddots & \vdots \\
a_n b_1 & a_n b_2 & a_n b_3    & \cdots & x_n \\
\end{bmatrix} .
\end{equation}
The correct expression for this determinant is
\begin{equation}
   \label{eq:1}
   P308 = \Big(\prod_{k=1}^n (x_k - a_k b_k)\Big) \Big(1 + \sum_{k=1}^n\frac{a_{k}b_{k}}{x_k - a_k b_k}\Big),
 \end{equation}
 which in Proskuryakov's book appears with denominators $x_{k}$
 instead of $x_{k}-a_{k}b_{k}$, see for instance \cite[page 321]{MR797072}.

 Indeed, this formula is a consequence of the more
 general one for the determinant of a sum of matrices \cite[pages
 162-163]{Marcus:fdma}, as it is also hinted in \cite[pages 40-41]{MR797072}. For convenience, we give here a self-contained
 proof based on the multilinearity of the determinant function.

 \begin{proof}[Proof of Formula~\eqref{eq:1}]
   Denote by $M$ the $n\times n$ matrix in \eqref{eq:4}, which  can be
   written as the sum of a diagonal and a rank 1 matrix as
   \begin{displaymath}
     M=\diag(x_{1}-a_{1}b_{1},\dots, x_{n}-a_{n}b_{n})+ a\cdot b^{T}
   \end{displaymath}
   for the $n$ vectors $a=(a_{1},\dots, a_{n})$ and
   $ b=(b_{1},\dots, b_{n})$.  Considering the determinant as a
   function of the columns of the matrix, we have that
   \begin{displaymath}
P308=     \det(M)= \det((x_{1}-a_{1}b_{1})\, e_{1} + b_{1} a,\dots, (x_{n}-a_{n}b_{n})\, e_{n}+b_{n} a),
   \end{displaymath}
   where $e_{i}$ denotes the standard $n$ vector
   $(0,\dots, 0, \overset{i}{1},0,\dots, 0)$. By the multilinearity of
   the determinant function and the fact that it vanishes when the
   vectors are linearly dependent, we have that
   \begin{align*}
          P308=&
                    \det((x_{1}-a_{1}b_{1})\, e_{1} ,\dots, (x_{n}-a_{n}b_{n})\, e_{n})\\
                  &+\sum_{k=1}^{n}
                    \det( b_{1} a,\dots, b_{k-1}a, (x_{k}-a_{k}b_{k})\, e_{k}, b_{k+1}a, \dots, b_{n}a) \\
     =& \prod_{k=1}^n (x_k - a_k b_k) + \sum_{i=1}^{n}a_{k}b_{k} \prod_{l\ne k} (x_{l} - a_{l} b_{l}),
   \end{align*}
   which gives the intended formula
        \end{proof}

        As an application, we compute the Ricci form of the
        Fubini-Study metric on the $n$-dimensional complex projective
        space $\mathbb{P}^{n}$. In Riemannian geometry, this
        computation is usually done using the invariance of this
        metric with respect to the action of the unitary group as in
        \cite[\S13.3]{Moroianu:lkg}. By contrast, Formula
        \eqref{eq:1} allows to do it in a direct way.

        Let $Z_{0},\dots, Z_{n}$ be the homogeneous coordinates of
        this projective space and for each $k\in \{0,\dots, n\}$
        consider the open chart
        $U_{k}=(Z_{k}\ne 0) \simeq \mathbb{C}^{n}$ with coordinates
        $z_{1},\dots,z_{n}$.  The \emph{Fubini-Study form}
        $\omega_{\FS}$ is the K\"ahler form on $\mathbb{P}^{n}$ given
        in these coordinates by
  \begin{displaymath}
    \omega_{\FS} := \iu\partial\overline\partial \log(1 + \|z\|^{2})
  \end{displaymath}
  where $\partial, \overline \partial$ are the Dolbeault operators and
  $\|z\|=(|z_{1}|^{2}+\cdots+|z_{n}|^{2})^{1/2}$. The corresponding
  Hermitian matrix with respect to the frame
  $\frac{\partial}{\partial z_{i}}$, $i=1,\dots, n$ writes down~as
  \begin{align*}
    H=&
    \begin{bmatrix}
\displaystyle{      \frac{\partial^{2}}{\partial z_{i}\partial \overline z_{j}}  \log(1 + \|z\|^{2}) }
    \end{bmatrix}_{i,j} \\
    =&
    \frac{1}{(1+\|z\|^{2})}
    \begin{bmatrix}
      1+\|z\|^{2}-\overline z_{1} \, z_{1}& -\overline z_{1} \, z_{2} &  \cdots & -\overline z_{1} \, z_{n}\\
      -\overline z_{2} \, z_{1} &       1+\|z\|^{2} -\overline z_{2} \, z_{2} & \cdots &       -\overline z_{2} \, z_{n}\\
      \vdots & \vdots & \ddots & \vdots \\
            -\overline z_{n} \, z_{1} &        -\overline z_{n} \, z_{2} & \cdots &       1+\|z\|^{2}-\overline z_{n} \, z_{n}
    \end{bmatrix}
  \end{align*}
and by \cite[Formula (12.6)]{Moroianu:lkg}, the  associated Ricci form is then given  by
  \begin{displaymath}
  \rho_{\FS}:= -\iu\partial\overline\partial \log(\det(H)).
  \end{displaymath}
  Notice that $\det(H)$ is a special case of $P308$ with
  \begin{displaymath}
x_i = \frac{1 + \|z\|^{2}  - |z_i|^2}{(1 + \|z\|^{2})^{2}},\quad
a_i = \frac{-\overline z_i}{(1 + \|z\|^{2})^{2}},\quad
b_i = \frac{z_i}{(1 + \|z\|^{2})^{2}} \quad \text{ for } i=1,\dots, n.
  \end{displaymath}
  Now a straightforward application of Formula \eqref{eq:1} gives
  $\det(H)= (1+\|z\|^{2})^{-n-1}$. This implies that
  \begin{displaymath}
    \rho_{\FS} = -\iu\partial\overline\partial \log( (1+\|z\|^{2})^{-n-1}) = (r+1) \,\omega_{\FS},
  \end{displaymath}
  showing that the Fubini-Study metric is Einstein with $r+1$ as Einstein constant.

\bibliographystyle{amsalpha}
\bibliography{./bibliography.bib}

\end{document}